\documentclass{amsart}
\usepackage{amssymb}
\usepackage[all]{xy}
\newtheorem{theorem}{Theorem}[section]

\newtheorem{lemma}[theorem]{Lemma}

\newtheorem{corollary}[theorem]{Corollary}

\theoremstyle{definition}

\newtheorem{definition}{Definition}[section]
\theoremstyle{remark}
\newtheorem{notation}[definition]{Notation}

\DeclareMathOperator\ort{\stackrel{\perp}{\oplus}}

\newcommand{\In}{[\,1\,;\,n\,]}
\newcommand{\orth}[1]{\boldsymbol{\omega}(#1)}
\newcommand{\field}[1]{\mathbb{#1}}
\newcommand{\bbbr}{\field{R}}
\newcommand{\bbbn}{\field{N}}
\newcommand{\rank}{{\rm rank}}
\newcommand{\Aut}{{\rm Aut}}

\newcommand{\mat}[1]{\boldsymbol{#1}}

\newcommand{\card}[1]{\bigl|\,{#1}\,\bigr|}

\newcommand{\dist}[2]{{\rm d}(#1,#2)}
\newcommand{\distf}[2]{{\rm d}_f(#1,#2)}
\DeclareMathOperator\distb{\rm dist}
\newcommand{\proj}{\bigl(\mat{I}-\frac{\mat{J}}{n}\bigr)}

\newcommand{\sidesetx}[3]{%
}
\newlength{\xxxxa}
\newlength{\xxxxb}

\begin{document}
\title{Regular Embeddings of Multigraphs}
\author{Hubert de Fraysseix}
\address{Centre d'Analyse et de Math\'ematiques Sociales\\
CNRS, UMR 8557\\
54 Bd Raspail, 75006 Paris\\
France}
\email{hf@ehess.fr}
\author{Patrice Ossona de Mendez}
\address{Centre d'Analyse et de Math\'ematiques Sociales\\
CNRS, UMR 8557\\
54 Bd Raspail, 75006 Paris\\
France}
\email{pom@ehess.fr}
\begin{abstract}
We prove that the vertex set of any twin-free loopless multigraph $G$ has an
embedding into some point set $P$ of some Euclidean space $\bbbr^k$,
such that the automorphism group of $G$ is isomorphic to the isometry
group of $\bbbr^k$ globally preserving $P$.
\end{abstract}
\maketitle
\section{Introduction}
\label{sec:intro}
Using spectral analysis, Babai proved in 1978 that the abstract
automorphism group of any multigraph $G$
having $s$ distinct eigenvalues with respective multiplicities
$m_1,m_2,\dots,m_s$ is a subgroup of
$\orth{m_1}\oplus\orth{m_2}\oplus\dots\oplus\orth{m_s}$,
where $\orth{m}$ denotes the real orthogonal group of dimension $m$
\cite{Babai}.
As a consequence, if all the eigenvalues of $G$ are simple, the only
automorphisms of $G$ are involutions.

Some years before, Mani proved that every triconnected planar graph $G$ can be
realized as the $1$-skeleton of a convex polytope $P$ in $\bbbr^3$
such that all automorphisms of $G$ are induced by isometries of $P$
\cite{Mani}.
One non trivial consequence of this result is that the automorphism
group $\Aut(G)$ of a triconnected planar graph $G$ has a chain
of normal subgroups $\Aut(G)=G_0\triangleright
G_1\triangleright\dots\triangleright G_m=1$ where each quotient
$G_i/G_{i-1}$ is either cyclic, or isomorphic to a symmetric group or $A_5$.

The result of Mani may be expressed in a weaker form: any triconnected
planar graph has an embedding $f$ into $\bbbr^3$, such that
$\Aut(G)$ is the
group of isometries of $\bbbr^3$ globally preserving the point set
$P=f(V(G))$, that we shall denote by $\orth{3,P}$.

These two results suggest that a graph $G$ may possibly have some
{\em regular embedding}, that is some embedding
 $f:V(G)\rightarrow\bbbr^k$ such that $\Aut(G)$ is
isomorphic to the group $\orth{k,f(V(G))}$ of isometries of $\bbbr^k$
globally preserving $f(V(G))$, and that this group might be expressed
as a subgroup of a group sum relying on spectral considerations. We
shall prove that this is indeed the case for any twin-free loopless multigraph.
In the remaining of this paper, multigraphs are always assumed to be loopless.

This result will be proved using techniques similar to the one used in
a graph symmetry detection heuristic presented in GD'99 \cite{Taxi_SymmGD99}.

In section~\ref{sec:def} we recall several definitions and introduce
notations. In section~\ref{sec:red} we reduce the study of the 
automorphism groups of multigraphs to the case of {\em irreducible}
multigraphs, that is to multigraphs having no twin vertices.
In section~\ref{sec:ed} we relate embeddings of a multigraphs to 
metrics and {\em distance matrices} on its vertex set. We reduce in section~\ref{sec:re} 
the problem of finding regular
embeddings to the one of finding metrics on the vertex set of the
multigraph which define {\em Euclidean, reconstructing} and 
  {\em commuting} distance matrices.
We prove in section \ref{sec:dist} that Euclidean distance matrices may be
built from symmetric real matrices with $0$ on the diagonal, what we
call {\em predistance matrices}.
 The built distance matrix is
commuting/reconstructing if the predistance is.  We deduce in 
section \ref{sec:rp} that any commuting reconstructing predistance
matrix defines a regular embedding and give some examples of reconstructing
commuting predistance matrices.
In section \ref{sec:regular} we study the special case of regular
multigraphs and give a strengthened version of Babai's result
\cite{Babai} in this context. Section~\ref{sec:cncl} is devoted to
concluding remarks.

\section{Definitions and Notations}
\label{sec:def}

Let $G$ be a multigraph. For any $x,y\in V(G)$, 
the {\em multiplicity} $\mu(x,y)$ is the number of edges (possibly
$0$) having $x$ and $y$ as endpoints.  
An {\em automorphism} of $G$ is
a one-to-one mapping $g:V(G)\rightarrow V(G)$ such that
$\mu(g(x),g(y))=\mu(x,y)$ for any $x,y\in V(G)$. The automorphisms of
$G$ define the {\em automorphism group} $\Aut(G)$ of $G$.

In order to ease the matrix presentation, we shall assume that the
vertex set of a multigraph $G$ of order $n$ is $\{1,\dots,n\}$.
Then the {\em adjacency matrix} of $G$ is the symmetric matrix with entries
$\mat{A}_{i,j}=\mu(i,j)$ and any automorphism $g\in\Aut(G)$ may be
written as a permutation matrix $\mat{g}$, where
$\mat{g}_{i,j}=\delta_{j,g(i)}$ ($\delta$ is the usual Kronecker
symbol). By permissive abridgment, $\Aut(G)$ would as well denote the
group of these permutation matrices. Notice that $\Aut(G)$ may be then
described as the group of the permutation matrices commuting with
$\mat{A}$. Also, we will denote by $\mat{E_i}$ the $n$ column matrix
whose $j$th entry is $\delta_{i,j}$, $\mat{1}$ the $n$ column matrix filled
with $1$s and by $\mat{J}$ the $n\times n$ matrix $\mat{1}\mat{1}^{\rm\sc t}$.

The symmetric group acting on a finite set $P$ will be denoted by $\mathfrak{S}(P)$.
The real orthogonal group of dimension $n$ will be denoted by $\omega(n)$.
If $P$ is a set of points of $\bbbr^n$, the subgroup of $\omega(n)$
globally preserving $P$ will be denoted by $\omega(n,P)$. Notice that 
$\omega(n,P)<\mathfrak{S}(P)$.

\section{Reducing Graphs}
\label{sec:red}
\begin{definition}
Let $G$ be a multigraph. $G$ is {\em reducible} if there exists two
vertices $x_1,x_2$ such that $\mu(x_1,y)=\mu(x_2,y)$ for any $y\in
V(G)\setminus\{x_1,x_2\}$. The vertices $x_1$ and $x_2$ are said to be
{\em twins}. If a multigraph is twin-free, it is {\em irreducible}.
\end{definition}
\begin{lemma}
Let $X$ be a subset of vertices, any two elements of which are twins.
Then, $\mu(x,y)$ is constant for $x\neq y\in X$.
\end{lemma}
\begin{proof}
Let $X=\{x_1,x_2,\dots,x_k\}$ Then, for $1<i<j\leq k$:
$\mu(x_i,x_j)=\mu(x_1,x_j)$ as $x_1$ and $x_i$ are twins and
$\mu(x_1,x_j)=\mu(x_1,x_2)$ as $x_2$ and $x_j$ are twins.
\end{proof}
\begin{lemma}
Let $G$ be a multigraph. For any distinct vertices $x,y,z$ the following holds:
\begin{align*}
\text{If $x$ and $y$ are twins, that is:}\quad\forall v\in V(G)\setminus\{x,y\},\quad &\mu(x,v)=\mu(y,v)\\
\text{and $y$ and $z$ are twins, that is:}\quad\forall v\in V(G)\setminus\{y,z\},\quad &\mu(y,v)=\mu(z,v)\\
\text{then $x$ and $z$ are twins, that is:}\quad\forall v\in V(G)\setminus\{x,z\},\quad &\mu(x,v)=\mu(z,v).
\end{align*}
\end{lemma}
\begin{proof}
Let $v\in V(G)\setminus\{x,z\}$.

If $v\neq y$:
$\mu(x,v)=\mu(y,v)$ (as $x$ and $y$ are twins) and
$\mu(y,v)=\mu(z,v)$ (as $y$ and $z$ are twins) hence
$\mu(x,v)=\mu(z,v)$.

Otherwise $\mu(x,z)=\mu(x,y)$ (as $y$ and $z$ are twins)
and $\mu(x,z)=\mu(y,z)$ (as $x$ and $y$ are twins) hence
$\mu(x,v)=\mu(z,v)$.
\end{proof}

\begin{corollary}
The vertex set of any multigraph $G$ has a unique partition
$\mathcal P=\{V_1,\dotsc,V_k\}$, such that any $V_i$ is a maximal subset
of twin vertices of $G$. This partition is the {\em
  twin-decomposition} of $G$.
\end{corollary}

\begin{corollary}
Let $G$ be a multigraph, let $\mathcal P$ be its twin decomposition
and let $G/\mathcal P$ denote the quotient multigraph.
Then 
\begin{equation}
\Aut(G)=\bigoplus_{X\in\mathcal P}\mathfrak{S}(X)\oplus\Aut(G/\mathcal
P)
\end{equation}
Moreover, no subset $Y\subseteq V(G/\mathcal P)$ of cardinality at
least two is such that $\mathfrak{S}(Y)$ is a subgroup of $\Aut(G/\mathcal P)$.
\end{corollary}
\begin{proof}
Assume there exists a subset $Y\subseteq V(G/\mathcal P)$ of cardinality at
least two is such that $\mathfrak{S}(Y)$ is a subgroup of
$\Aut(G/\mathcal P)$. Let $\tau$ be a transposition in
$\mathfrak{S}(Y)$ exchanging two vertices $a$ and $b$ of $G/\mathcal
P$. As $\tau$ is an automorphism of $G/\mathcal P$ it follows that $a$
and $b$ are twins in $G/\mathcal P$. Identifying $a$ and $b$ with the
classes of twins of vertices in $G$, it follows that any $x\in a$ is a
twin of any $y\in b$ in $G$ hence $a\cup b$ is a class of twins, what
contradicts the maximality of the classes in $\mathcal P$.
\end{proof}

\section{Embeddings and Distances} 
\label{sec:ed}
Recall that a {\em metric} on a set $X$ is a mapping ${\rm d}:
X^2\rightarrow \bbbr^+$ satisfying the usual axioms of a metric,
that is:
\begin{align*}
\forall (x,y)\in X^2, \qquad &\dist{x}{y}=\dist{y}{x}\\
\forall (x,y)\in X^2, \qquad &\dist{x}{y}=0\iff x=y\\
\forall (x,y,z)\in X^3, \qquad &\dist{x}{z}\leq \dist{x}{y}+\dist{y}{z}
\end{align*}

An {\em Euclidean metric} on a set $X$ is a metric
${\rm d}$ on $X$ such that there exist
some Euclidean space $\bbbr^k$ and an embedding
$f:X\rightarrow\bbbr^k$ so that, for any $x,y\in X$:
\begin{equation*}
\distb(f(x),f(y))=\dist{x}{y}
\end{equation*}
where $\distb$ is the Euclidean metric of $\bbbr^k$. 
Any embedding of a multigraph $G$ into some Euclidean space
$\bbbr^k$ defines an Euclidean metric ${\rm d}_f$ on $V(G)$ by
$\distf{x}{y}=\distb(f(x),f(y))$. 

\begin{definition}
The {\em distance matrix} of a metric ${\rm d}$ on the set
$\{1,\dots,n\}$ is the $n\times n$ real symmetric matrix $\mat{D}$ defined by:
\begin{equation*}
\forall i,j\in\{1,\dots,n\},\qquad\mat{D}_{i,j}=\dist{i}{j}^2
\end{equation*}
\end{definition}
Notice that the entries of $\mat{D}$ are the {\bf squares} of the
distances and not the distances themselves. 

\begin{definition}
Let $G$ be a multigraph  
and let $\mat{D}$ be a distance matrix of a metric defined on $V(G)=\{1,\dots,n\}$.

The distance matrix $\mat{D}$ is
\begin{itemize}
\item {\em Euclidean} if the metric from which $\mat{D}$ comes from is
  Euclidean. A compatible embedding of $G$ into an Euclidean space
  $\bbbr^k$ is then called a {\em $\mat{D}$-embedding} of $G$;
\item {\em reconstructing} if $G$ may be reconstructed from $\mat{D}$, that is:
\begin{equation*}
\exists\Xi:\bbbr\rightarrow\bbbn,\quad \forall i\neq j\in
\{1,\dots,n\},\quad\mu(i,j)=\Xi(\mat{D}_{i,j});
\end{equation*}
\item {\em commuting} if any automorphism of $G$ commutes with
  $\mat{D}$, that is:
\begin{align*}
\forall g\in\Aut(G),\qquad &\mat{g}\mat{D}=\mat{D}\mat{g}
\intertext{or, equivalently:}
\forall g\in\Aut(G),\forall
i,j\in\{1,\dots,n\}\qquad &\mat{D}_{g(i),g(j)}=\mat{D}_{i,j}.
\end{align*}
\end{itemize}
\end{definition}
\section{Distance Matrices and Regular Embeddings}
\label{sec:re}

\begin{definition}
Let $G$ be a multigraph and let $k$ be an integer.
A {\em regular embedding} of $G$ into $\bbbr^k$ is a mapping
$f:V(G)\rightarrow \bbbr^k$ such that 
$\Aut(G)$ is isomorphic to $\omega(k,f(V(G)))$.
\end{definition}

We recall the classical following theorem (which proof may be found in
any undergraduate textbook):
\begin{theorem}[Isometry Extension Theorem]
\label{th:iet}
Let $(a_i)_{i\in I}$ and $(b_i)_{i\in I}$ be families of points of an
Euclidean affine space $\mathcal E$ (with distance $\distb$) such that:
\begin{align*}
\forall (i,j)\in I^2\qquad&\distb(a_i,a_j)=\distb(b_i,b_j)
\intertext{Then there exists an isometry $f$ of $\mathcal E$ such
  that:}
\forall i\in I\qquad&f(a_i)=b_i
\end{align*}
Moreover, if $(b_i)_{i\in I}$ spans $\mathcal E$ then $f$ is uniquely determined.
\end{theorem}

\begin{theorem}
\label{thm:reg}
Let $f$ be a one-to-one embedding of a multigraph $G$ into the
Euclidean space $\bbbr^k$ (with metric $\distb$) 
and let $\mat{D}$ be the Euclidean distance matrix of the metric
defined by $f$ on $V(G)$. 

Assume $\mat{D}$ is reconstructing.
Then $f$ is regular if and only if $f(V(G))$ spans $\bbbr^k$ and $\mat{D}$
is commuting. 
\end{theorem}
\begin{proof}
Assume $f(V(G))$ spans $\bbbr^k$ and that $\mat{D}$ is commuting.

Let $\distb$ be the usual Euclidean metric of $\bbbr^k$.
Let $g\in\Aut(G)$. Then, for any vertices $x,y$ of $G$,
$\distb(f(g(x)),f(g(y)))=\distb(f(x),f(y))$.
According to Theorem \ref{th:iet},
$g$ may be extended to a unique
isometry of $\bbbr^k$ preserving $f(V(G))$,
that is: to a unique element of $\omega(k,f(V(G))$. 

Let $\tilde{g}$ denote the isometry defined by the automorphism $g$.
Obviously, if $g_1,g_2\in\Aut(G)$,
$\widetilde{g_1g_2}=\tilde{g_1}\tilde{g_2}$. Thus we have defined a
group morphism from $\Aut(G)$ to $\omega(k,f(V(G)))$.

Now assume $\phi\in\omega(k,f(V(G)))$.
Define $g:V(G)\rightarrow V(G)$ by
$f(g(x))=\phi(f(x))$. Notice that $g$ is well defined as $\phi$ is
obviously one-to-one.
As $\mat{D}$ is reconstructing, there exists a mapping
$\Xi:\bbbr\rightarrow\bbbn$ such that $\Xi(\mat{D}_{i,j})=\mu(i,j)$.
 Let $i\neq j\in\{1,\dots,n\}$. 
Then, for any $i\neq j\in\{1,\dots,n\}$, we have:
$\mu(i,j)=\Xi(\mat{D}_{i,j})$. As $\phi$ is an isometry, the distance
between $f(i)$ and $f(i)$ is the same as the distance as the distance
between $\phi(f(i))$ and $\phi(f(i))$, that is: between $f(g(i))$ and
$f(g(i))$ (by the definition of $g$). It follows that
$\mat{D}_{g(i),g(j)}=\mat{D}_{i,j}$ hence
$\mu(g(i),g(j))=\mu(i,j)$. It follows that
$g$ is an automorphism of $G$. Moreover, $\phi$ extends
$f^{-1}\circ g\circ f$, thus $\phi=\tilde{g}$. It follows that
$g\mapsto \tilde{g}$ is actually a group isomorphism from $\Aut(G)$ to
$\omega(k,f(V(G)))$.

Conversely, assume $f$ is a regular embedding. Then $f(V(G))$ spans $\bbbr^k$ for
otherwise $\omega(k,f(V(G)))$ would not be a finite group although
$\Aut(G)$ is. Let $\phi\in\omega(k,f(V(G)))$.
Define $\tilde{\phi}:V(G)\rightarrow V(G)$ by
$f(\tilde{\phi}(x))=\phi(f(x))$. As previously, 
we get $\mu(\tilde{\phi}(x),\tilde{\phi}(y))=m$, hence
$\tilde{\phi}$ is an automorphism of $G$ such that
$\dist{f(\tilde{\phi}(x))}{f(\tilde{\phi}(y))}=\dist{f(x)}{f(y)}$.
Moreover, $\widetilde{\phi_1\phi_2}=\tilde{\phi_1}\tilde{\phi_2}$. It
follows that $\phi\mapsto\tilde{\phi}$ is a group-morphism. As it is
clearly one-to-one and as $\Aut(G)$ and $\omega(k,f(V(G)))$ are
isomorphic by assumption, $\phi\mapsto\tilde{\phi}$ is onto. It
follows that $\dist{f(g(x))}{f(g(y))}=\dist{f(x)}{f(y)}$ for any
$g\in\Aut(G)$, that is: $\mat{D}$ is commuting.
\end{proof}

We shall consider now this result from another point of view:

\begin{corollary}
Let $G$ be a multigraph and let $\mat{D}$ be an Euclidean
reconstructing and commuting distance matrix on $G$. Then $\mat{D}$
defines a regular embedding of $G$.
\end{corollary}
\begin{proof}
Consider any $\mat{D}$-embedding $f$ of $G$ and the subspace spanned
by $f(V(G))$ and apply Theorem~\ref{thm:reg}.
\end{proof}

Our main problem is now to build such a distance matrix.
\section{Euclidean distance matrices from Predistance matrices}
\label{sec:dist}
Given a distance $\distb$ on the set $\{1,\dots,n\}$, it is classical to define the
corresponding bilinear form $\mat{B}$ by:
\begin{equation}
\mat{B}=-\frac{1}{2}(\mat{I}-\frac{1}{n}\mat{J})\mat{D}(\mat{I}-\frac{1}{n}\mat{J})
\end{equation}
where $\mat{D}_{i,j}=\distb(i,j)^2$.

It is well known that the distance $\distb$ is Euclidean (i.e. allows
an isometric embedding into some Euclidean space) if and only if
$\mat{B}$ is positive semi-definite \cite{Benzecri2,Hakimi}.
Such a result extends to a characterization of those symmetric real
matrices which are Euclidean distance matrices.

\begin{definition}
A {\em predistance matrix} on $G$ is an $n\times n$ symmetric real matrix $\mat{P}$
with $0$ on the diagonal (notice that negative entries are allowed).

The {\em bilinear form} $\Lambda(\mat{P})$ {\em of} $\mat{P}$ is
defined by:
\begin{equation*}
\Lambda(\mat{P})=-\frac{1}{2}(\mat{I}-\frac{1}{n}\mat{J})\mat{P}(\mat{I}-\frac{1}{n}\mat{J})
\end{equation*}
\end{definition}

\begin{lemma}
\label{lem:B}
Let $\mat{D}$ be a predistance matrix and let $k$ be an integer.
Then, the following statements are equivalent:
\begin{enumerate}
\item $\mat{D}$ is an Euclidean distance matrix on $G$,
and there exists a $D$-embedding of $G$ in $\bbbr^k$;
\item the matrix $\Lambda(\mat{D})$
 is positive semi-definite and $\rank(\Lambda(\mat{D}))\leq k$.
\end{enumerate}
\end{lemma}
\begin{proof}
Assume $\Lambda(\mat{D})$ is positive semi-definite. Then the square of the distance 
defined by $\Lambda(\mat{D})$ between the basis points $e_i$ and $e_j$ is
$(\mat{E}_i-\mat{E}_j)^{\sc t}\Lambda(\mat{D})(\mat{E}_i-\mat{E}_j)$. As
$\proj(\mat{E}_i-\mat{E}_j)=(\mat{E}_i-\mat{E}_j)$, we get
that the square of this distance equals
$(\mat{D}_{i,j}+\mat{D}_{j,i}-\mat{D}_{i,i}-\mat{D}_{j,j})/2=\mat{D}_{i,j}$
as $\mat{D}$ is symmetric with zero on the diagonal. We define $f(i)$
as the projection of $e_i$ orthogonal to the isotropic space of
$B$.

Conversely, assume $\mat{D}$ is an Euclidean distance matrix and let
$f:\In\rightarrow\bbbr^k$ be a $\mat{D}$-embedding of $G$.
Then $\Lambda(\mat{D})$ is the Gram matrix of the vectors $f(i)$,
i.e. $\Lambda(\mat{D})_{i,j}=<f(i),f(j)>$.
It is well known that this Gram matrix determines the vectors $f(i)$ up
to isometry.
\end{proof}

Let $\mat{P}$ be a predistance matrix.
The condition that $\Lambda(\mat{P})$ has no negative eigenvalue is
quite difficult to handle. However, we know the following about $\Lambda(\mat{P})$:
\begin{itemize}
\item $\Lambda(\mat{P})$ is symmetric real and thus diagonalizable,
\item $\mat{1}\in\ker(\Lambda(\mat{P}))$ as $(\mat{I}-\frac{1}{n}\mat{J})\mat{1}=\mat{0}$,
\item $\Lambda(\mat{P})$ has an orthogonal basis $\mathcal B(\Lambda(\mat{P}))$ of eigenvectors
  including $\mat{1}$. The eigenvectors of $\mathcal
  B(\Lambda(\mat{P}))\setminus\{\mat{1}\}$ have eigenvalues
  $\lambda_1>\dots>\lambda_r$ with respective multiplicities
  $m_1,\dots,m_r$.
\end{itemize}

Now consider the {\em reduced distance matrix}
$\mat{P^\star}=\mat{P}-2\lambda_r(\mat{J}-\mat{I})$.

\begin{lemma}
\label{lem:Bstar}
 $\mathcal{B}(\Lambda(\mat{P}))$ is an orthogonal basis of
eigenvectors of $\Lambda(\mat{P^\star})$, and the eigenvalue (for
$\Lambda(\mat{P^\star})$) of 
$v\in\mathcal{B}(\Lambda(\mat{P}))$ is 
$$\begin{cases}
0&\text{if }\mat{v}=\mat{1}\\
\lambda-\lambda_r&\text{if }\mat{v}\neq\mat{1} 
\text{ and } \Lambda(\mat{P})\mat{v}=\lambda \mat{v}
\end{cases}$$
Thus, $\Lambda(\mat{P^\star})$ is positive semi-definite and has corank $m_r+1$.
\end{lemma}
\begin{proof}
If $\mat{v}\neq\mat{1}$, then $\mat{v}\in\mat{1}^\perp$, and thus
$\proj \mat{v}=\mat{v}$.
Hence, if $\Lambda(\mat{P})\mat{v}=\lambda \mat{v}$, then
\begin{align*}
\Lambda(\mat{P^\star})\mat{v}&=-\frac{1}{2}\proj\mat{P^\star}\proj \mat{v}\\
&=-\frac{1}{2}\proj(\mat{P}-2\lambda_r(\mat{J}-\mat{I}))\proj \mat{v}\\
&=-\frac{1}{2}\proj\mat{P}\proj
\mat{v}+\lambda_r\proj(\mat{J}-\mat{I})\proj\mat{v}\\
&=\Lambda(\mat{P})\mat{v}-\lambda_r\mat{v}&\text{(as
}\mat{J}\mat{v}=\mat{0}\text{)}\\
&=(\lambda-\lambda_r)\mat{v}
\end{align*}
\end{proof}

\section{Reconstructing and Commuting predistance matrices}
\label{sec:rp}
\begin{definition}
A predistance matrix $\mat{P}$ is
\begin{itemize}
\item {\em commuting} if $\mat{P}$ commutes with any automorphism
  of $G$, that is:
\begin{align*}
\forall g\in\Aut(G),\qquad&\mat{g}\mat{P}=\mat{P}\mat{g}
\intertext{or, equivalently:}
\forall g\in\Aut(G), \forall
i,j\in\{1,\dots,n\},\qquad&\mat{P}_{g(i),g(j)}=\mat{P}_{i,j};
\end{align*}
\item {\em reconstructing} if there exists a mapping
  $\Xi:\bbbr\rightarrow\bbbn$ such that $\mu(i,j)=\Xi(\mat{P}_{i,j})$,
  or equivalently:
\begin{equation*}
\forall i,j,i',j'\in\{1,\dots,n\}\qquad \mat{P}_{i,j}=\mat{P}_{i',j'}\Rightarrow\mu(i,j)=\mu(i',j')
\end{equation*}
\end{itemize}
\end{definition}

\begin{notation}
Let $\mat{P}$ be a predistance matrix. As $\Lambda(\mat{P})$ has $0$ as an eigenvalue, we may define
$\zeta(\mat{P})$ as follows:
\begin{itemize}
\item if the smallest eigenvalue $\lambda$ of $\Lambda(\mat{P})$ is negative
  and has multiplicity $m$ then $\zeta(\mat{P})=n-m-1$,
\item if the smallest eigenvalue of $\Lambda(\mat{P})$ is $0$ with multiplicity
  $m>1$ then  $\zeta(\mat{P})=n-m$,
\item if the smallest eigenvalue of $\Lambda(\mat{P})$ is $0$ with multiplicity
  $1$ and if the second smallest eigenvalue of $\Lambda(\mat{P})$ has
  multiplicity $m$ then $\zeta(\mat{P})=n-m-1$.
\end{itemize}
\end{notation}

\begin{theorem}
\label{thm:p2r}
Let $G$ be an irreducible multigraph.
Any commuting reconstructing predistance matrix $\mat{P}$ defines a regular
embedding $f$ into $\bbbr^{\zeta(\mat{P})}$. 
\end{theorem}
\begin{proof}
Let $\mat{P}$ be a commuting reconstructing predistance matrix $\mat{P}$. 
The corresponding reduced distance matrix
$\mat{P^\star}$, according to Lemma~\ref{lem:B} and
Lemma~\ref{lem:Bstar} is Euclidean and there exists a $\mat{P^\star}$
embedding $f$ on $G$ into $\bbbr^{\zeta(\mat{P})}$.  
As $\mat{P}$ is reconstructing and commuting, so is $\mat{P^\star}$.
Assume $f(x)=f(y)$. Since $f$ is reconstructing, $x$ and $y$ are twins
of $G$, contradicting its irreducibility. Hence $f$ is one-to-one and,
according to Theorem~\ref{thm:reg}, $f$ is a regular embedding.
\end{proof}
\begin{corollary}
Let $G$ be an irreducible multigraph and let $\mat{P}$ be a
reconstructing predistance matrix which commutes with the automorphisms of
$G$.

Denote by $\mat{1}\bbbr$ the line spanned by $\mat{1}$ and by $\ort$
the orthogonal direct sum of vector spaces.
Let $\mat{1}\bbbr\ort E_1\ort E_2\ort\dots\ort E_r$ be a decomposition
into eigenspaces of $\Lambda(\mat{P})$, the eigenvalues associated with
$E_1,\dots,E_r$ being $\lambda_1>\lambda_2>\dots>\lambda_r$ (notice
that some $\lambda_i$ may be $0$).
Then the abstract automorphism group of $G$ is a subgroup of
$$\orth{\dim E_1,P_1}\oplus\orth{\dim
  E_2,P_2}\oplus\cdots\oplus\orth{\dim E_{r-1},P_{r-1}}$$
where $P_i$ is the orthogonal projection of the image of $G$ under the
regular embedding defined by $\mat{P}$ into the subspace $E_i$.
\end{corollary}

Here some examples of simple commuting reconstructing predistance matrices
of a simple connected graph $G$ (the two last examples only apply if
$G$ has order at least $3$):

\begin{align*}
  \text{Adjacency}\qquad&\mat{A}_{i,j}=\begin{cases}
    0,\qquad&\text{if $i=j$ or $i$ and $j$ are adjacent}\\
    1,&\text{otherwise}
  \end{cases}\\
  \text{Graph distance}\qquad&\mat{S}_{i,j}=\begin{cases}
    0,\qquad&\text{if $i=j$}\\
    l,&\text{if $l$ is the graph distance from $i$ to $j$}
  \end{cases}\\
  \text{Bisected Czekanovski-Dice}\qquad&\mat{C}_{i,j}=\begin{cases}
    0,&\text{if $i=j$}\\
    1-\frac{2}{d(i)+d(j)},&\text{if $i$ and $j$ are adjacent}\\
    1,&\text{otherwise}
  \end{cases}\\
  \text{Q-distance}\qquad&\mat{Q}_{i,j}=\begin{cases}
    0,&\text{if $i=j$}\\
    1,&\text{if $i$ and $j$ are non adjacent}\\
    1-\frac{1}{\sqrt{d(i)d(j)}},&\text{otherwise}
  \end{cases}\\
\end{align*}

\begin{theorem}
Any irreducible multigraph has a regular embedding into some Euclidean
space.
\end{theorem}
\begin{proof}
The adjacency matrix $\mat{A}$ of $G$ defined by
$\mat{A}_{i,j}=\mu(i,j)$ is commuting and reconstructing. Hence the result
follows from Theorem~\ref{thm:p2r}.
\end{proof}

\section{Regular Multigraphs}
\label{sec:regular}
As noted in section~\ref{sec:rp}, the adjacency matrix $\mat{A}$
defines a commuting reconstructing predistance matrix.
In the particular case where $G$ is $d$-regular, we have
$\mat{A}\mat{J}=d\mat{J}$ thus $\Lambda(\mat{A})=(1-d/n)^2\mat{A}$. 
Moreover (see \cite{agt}):
\begin{itemize}
\item $d$ is an eigenvalue of $G$ with eigenvector $\mat{1}$.
\item If $G$ is connected, then the multiplicity of $d$ is one.
\item For any eigenvalue $\lambda$ of $G$, we have $\card{\lambda}\leq
  d$.
\end{itemize}
In this case we deduce a strengthening of Babai's result \cite{Babai}
mentioned in section~\ref{sec:intro}:
\begin{corollary}
Let $G$ be a connected irreducible regular multigraph having $s$
distinct eigenvalues $\lambda_1>\lambda_2>\dots>\lambda_s$ with
respective multiplicities $m_1,m_2,\dots,m_s$. Then, the abstract
automorphism group $\gamma$ of $G$ is a subgroup of
$$\orth{m_2,P_2}\oplus\orth{m_3,P_3}\oplus\dots\oplus\orth{m_{s-1},P_{s-1}}.$$
\end{corollary}
\begin{proof}
Assume $G$ is connected. Let $\mathcal B(\mat{A})$ be an orthogonal
basis of eigenvectors of $\mat{A}$. As $\mat{J}\mat{v}=0$ for any
vector orthogonal to $\mat{1}$ and as  $\mat{J}\mat{1}=n\mat{1}$,
$\mathcal B(\mat{A})$ is also a basis of eigenvectors of $\Lambda(\mat{A})$.
The eigenvalue (for $\Lambda(\mat{A})$) of an eigenvector $\mat{v}$ with
eigenvalue $\lambda$ (for $\mat{A}$) is clearly $\lambda$ if
$\lambda\neq d$ and $0$ if $\lambda=d$. As the eigenvalue $d$ is
maximal and simple for $\mat{A}$, we deduce that
$\zeta(\mat{A})=n-m-1$, where $m$ is the multiplicity of the smallest
eigenvalue of $\mat{A}$. Moreover, the greatest eigenvalue is simple
and corresponds to eigenvector $\mat{1}$, according to
Perron-Frobenius theorem \cite{frob,perron}.
\end{proof}

As a matter of fact, 
instead of
$\orth{m_2,P_2}\oplus\orth{m_3,P_3}\oplus\dots\oplus\orth{m_{s-1},P_{s-1}}$,
we may prove a similar result with
$\orth{m_3,P_3}\oplus\dots\oplus\orth{m_{s-1},P_{s-1}}\oplus\orth{m_{s},P_{s}}$
by considering
$2\lambda_2(\mat{J}-\mat{I})-\mat{A}$ instead of
$\mat{A}-2\lambda_s(\mat{J}-\mat{I})$.

A {\em strongly regular graph} is a regular graph such that there
exist constants $\lambda,\mu$ so that any two adjacent vertices have
exactly $\lambda$ common neighbors, and every two non-adjacent
vertices have exactly $\mu$ common neighbors \cite{seidel}.
As a strongly regular graph has exactly three distinct eigenvalues \cite{agt}, we
may choose to keep among $E_2$ and $E_3$ the one having the smallest dimension.
Hence every connected irreducible strong regular graph
of order $n$ has a regular embedding into an Euclidean space of
dimension at most $\lfloor\frac{n-1}{2}\rfloor$. As an example, the
Petersen graph has a regular embedding into $\bbbr^4$, which is
optimal (the automorphism group of the Petersen graph cannot be
realized as the isometry group of a set of points in $\bbbr^3$).

\section{Conclusion}
\label{sec:cncl}
The techniques presented here allow to construct regular embeddings
from weakly constrained matrices: we only ask that they should be
symmetric, have $0$ on the diagonal, are reconstructing and commute
with every automorphism. 

Recall that a {\em cellular algebra}, or {\em coherent algebra}, is an
algebra of $n\times n$ complex matrices which has a basis
$\{\mat{B}_0,\mat{B}_1,\dots,\mat{B}_t\}$ consisting of matrices with
entries $0$ and $1$ satisfying the following conditions:
\begin{enumerate}
\item $\mat{B}_0+\mat{B}_1+\cdots+\mat{B}_t=\mat{J}$;
\item $\exists 0\leq r\leq t:\quad \mat{B}_0+\mat{B}_1+\cdots+\mat{B}_r=\mat{I}$;
\item the set $\{\mat{B}_0,\mat{B}_1,\dots,\mat{B}_t\}$ is closed
  under transposition.
\end{enumerate}
The unique minimal cellular algebra which contains $\mat{A}$ as an element is
the {\em cellular algebra generated by} $\mat{A}$ and is denoted by
$<<\mat{A}>>$. Notice that a basis of this cellular algebra may be
constructed in polynomial time. In the case where $\mat{A}$ is the
adjacency matrix of a simple graph $G$, the cellular algebra
$<<\mat{A}>>$ is called the cellular algebra {\em generated by} $G$,
or the cellular algebra {\em of} $G$. This algebra may contain
non symmetric basis elements. However, any matrix in $<<A>>$
commutes with any automorphism of $G$. It follows that
 any commuting predistance matrix $\mat{P}$ in
$<<A>>$ may be written as $\sum_{i=r+1}^t
\lambda_{i}(\mat{B}_i+\mat{B}^{\rm\sc t}_i)$, with $\lambda_i\in\bbbr$
for $r< i\leq t$ and thus form an affine space $\mathcal E$.
By perturbing each of them (by adding $\epsilon(\mat{P})\mat{A}$, for instance) we
obtain reconstructing and commuting predistance matrices forming a
dense subset of $\mathcal E$.

We hope that the suggested approach will be fruitful and more
practical than the usual techniques arising from algebra and spectral
analysis.

\providecommand{\bysame}{\leavevmode\hbox to3em{\hrulefill}\thinspace}
\providecommand{\MR}{\relax\ifhmode\unskip\space\fi MR }
\providecommand{\MRhref}[2]{%
  \href{http://www.ams.org/mathscinet-getitem?mr=#1}{#2}
}
\providecommand{\href}[2]{#2}

\end{document}